\documentclass[11pt,a4paper]{article} 

\usepackage{amsmath} 
\usepackage[T1]{fontenc} 
\usepackage[utf8]{inputenc} 

\usepackage[ngerman,french,italian]{babel}
\usepackage{geometry} 
\usepackage{parskip}
\usepackage{titling}
\usepackage{sectsty}
\usepackage{graphicx}
\usepackage{float}
\usepackage{hyperref}
\usepackage{caption}

\usepackage{xparse}

\NewDocumentCommand{\INTERVALINNARDS}{ m m }{
    #1 {,} #2
}
\NewDocumentCommand{\interval}{ s m >{\SplitArgument{1}{,}}m m o }{
    \IfBooleanTF{#1}{
        \left#2 \INTERVALINNARDS #3 \right#4
    }{
        \IfValueTF{#5}{
            #5{#2} \INTERVALINNARDS #3 #5{#4}
        }{
            #2 \INTERVALINNARDS #3 #4
        }
    }
}

\usepackage{xcolor} 

\usepackage{enumerate} 

\usepackage{cancel}

\addto\captionsngerman{} 
\pretitle{\raggedright\LARGE} 
\posttitle{\par\vskip 0.5em} 

\preauthor{\raggedright\large} 
\postauthor{\par\vskip 0.5em}

\allsectionsfont{\bfseries\sffamily} 

\title{{\sffamily\bfseries Una dimostrazione diretta della legge di probabilità di Poisson (A direct proof of the Poisson probability law)}} 
\author{{\sffamily Pier Franco Nali, Cagliari (IT), pfnali@alice.it}} 

\date{} 

\makeatletter 
\newcommand{\xRightarrow}[2][]{\ext@arrow 0359\Rightarrowfill@{#1}{#2}}   
\makeatother 

\begin{document}

\selectlanguage{italian} 

\clearpage\maketitle 
\thispagestyle{empty}

The purpose of this paper is to prove directly, by an elementary method, the Poisson probability law. This proof is offered as an alternative to the more usual derivation from binomial distribution in the limit of small probabilities. The same proof is then applied to the solution of a problem in statistical mechanics.

Lo scopo di questo articolo è dimostrare direttamente, con un metodo elementare, la legge di probabilità di Poisson. Questa dimostrazione è proposta in alternativa alla più consueta derivazione dalla distribuzione binomiale nel limite delle piccole probabilità. La stessa dimostrazione viene quindi applicata alla soluzione di un problema di meccanica statistica.
\section{Introduzione}
La legge (o funzione, o distribuzione\footnote{C'è una certa ambiguità dei termini: i fisici li usano talvolta come sinonimi, i matematici preferiscono parlare di funzione di distribuzione (o ripartizione) quando si riferiscono alle probabilità cumulate. Qui seguiremo l'uso dei fisici.}) di probabilità di Poisson, o poissoniana,
\[P(x;\mu)=\frac{\mu\,^xe^{-\mu}}{x\,!}\]
di parametro \(\mu\,\), esprime la probabilità che \(x\) eventi si verifichino in maniera casuale in un dato intervallo di tempo quando in media se ne verificano \(\mu\,\). La si utilizza in situazioni in cui degli eventi si susseguono nel tempo in modo indipendente come, ad es., le disintegrazioni radioattive o le chiamate entranti in un centralino telefonico. È stato osservato \cite{hodges1960} che molti libri di testo di teoria delle probabilità si accontentano di derivare la poissoniana come limite della distribuzione binomiale
\[B(x;n,p)=\binom{n}{x}p\,^x q\,^{n-x}\,,\]
dove \(x\) è in numero di successi su \(n\) prove con probabilità di successo e insuccesso, rispettivamente, \(p\) e \(q\,\). La derivazione procede più o meno in questo modo: C'è un gran numero \(n\) di eventi che possono accadere (ad es., molti abbonati al telefono possono fare una chiamata in un certo intervallo di tempo) e c'è una piccola probabilità \(p\) che uno specifico evento tra questi accada (ad es., che uno specifico abbonato chiami). Assumendo che gli eventi siano indipendenti, il numero di eventi (ad es., le chiamate) ha esattamente la distribuzione binomiale \(B(x;n,p)\,\). Facendo tendere \(n\) all'infinito e \(p\) a zero in modo tale da mantenere costante il prodotto \(n\times p=\mu>0\,\), viene mostrato che la binomiale tende alla poissoniana \(P(x;\mu)\) di valore atteso \(\mu\,\). Questo risultato, chiamato teorema (o limite) di Poisson, è molto importante e molto noto. È in questo modo che storicamente è stata derivata la poissoniana \cite{poisson1837} ed è stata utilizzata nelle prime applicazioni pratiche \cite{bortkiewicz1898}. In seguito si è compreso che il terorema di Poisson è valido, in una versione modificata, anche in situazioni in cui gli eventi non hanno tutti la stessa probabilità, e la distribuzione ha una forma più generale di quella originale (distribuzione Poisson--binomiale) \cite{mises1921}. 

\begin{figure}[H]
\includegraphics[width=\textwidth]{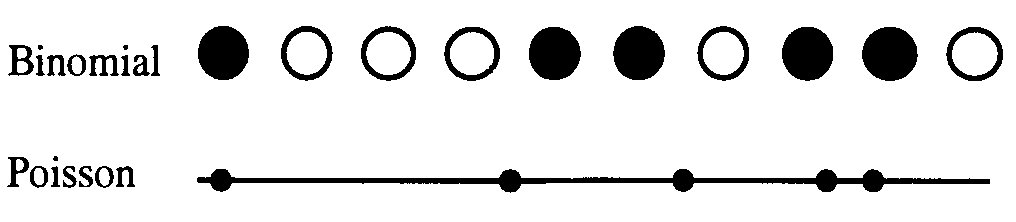}
\caption{Raffronto schematico tra binomiale e poissoniana \cite{hutchinson1997}}\label{fig:PointProcess}
\centering
\end{figure}

La stretta relazione tra binomiale e poissoniana mette in risalto la similarità delle due distribuzioni. Vi è però tra di esse un'importante differenza, per non dire un contrasto \cite{hutchinson1997}: con la binomiale, l'uno o l'altro di due tipi di eventi (ad es., successo o insuccesso) si verificano in particolari occorrenze di una sequenza discreta (ad es., in una serie di tentativi); con la poissoniana, l'evento rilevante avviene in punti casuali di un continuo (fig. \ref{fig:PointProcess}). 
Un'altra evidenza che si tratta di due distribuzioni distinte è di tipo storico \cite{duker1955}: dopo Poisson, la distribuzione che porta il suo nome è stata "riscoperta" da numerosi altri investigatori utilizzando svariati metodi. Essa può essere derivata o come limite della binomiale o come distribuzione indipendente di eventi rari distribuiti a caso. Le derivazioni indipendenti sono riportate meno di frequente nei libri di testo e generalmente presuppongono conoscenze più avanzate rispetto alla consueta derivazione dalla binomiale.

Lo scopo di questo articolo è presentare una derivazione indipendente della poissoniana, ottenuta direttamente - e in modo elementare - dal processo di Poisson\footnote{Con questo termine si indica un particolare processo stocastico che fornisce, con opportune varianti, il modello matematico di moltissimi fenomeni fisici, sociali e demografici caratterizzati da eventi che si susseguono nel tempo in modo aleatorio.} in una dimensione (linea temporale). Questa viene poi generalizzata a processi in più dimensioni e applicata alla soluzione di un classico problema di meccanica statistica (distribuzione spaziale delle molecole di un gas).
\section{Processo di Poisson e poissoniana \label{sez_2}}
Un processo di Poisson è caratterizzato, nel caso più semplice, da una successione di eventi indipendenti con tasso \(\lambda\) costante (in media) su un intervallo continuo di estensione \(t>0\) (periodo di osservazione). In situazioni più complesse, il flusso degli eventi può presentare variazioni rispetto al semplice schema stazionario (ad es., il tasso può variare in funzione del tempo). Nella rappresentazione grafica molto intuitiva della fig. \ref{fig:PointProcess} gli eventi in successione sono raffigurati come punti isolati su una linea, che può simboleggiare un range temporale o un tratto della retta dei numeri reali positivi. Questa descrizione unidimensionale può essere generalizzata a due, tre o più dimensioni considerando punti esistenti in  aree, volumi o insiemi misurabili in genere \cite{wiki:xxx}. Tornando al processo poissoniano del tipo più semplice, detto "stazionario a tasso medio costante", vale la condizione
\begin{equation}\label{eq_1}
\lim_{t\to\infty}\frac{x}{t}=\lambda\,,    
\end{equation}
dove \(x\) è il numero cumulato di eventi conteggiati nel tempo \(t\,\).

Ma come si può stabilire se un determinato processo ha, nel concreto, un tasso stazionario alla scala temporale accessibile all'osservazione o all'esperimento? Il solo modo è ripetere il conteggio del numero di eventi \(x_i\) in intervalli (non necessariamente equiestesi) di estensione \(t_i\) e determinare se c'è una tendenza stazionaria nella serie dei rapporti \(x_i/t_i\). Poiché questi rapporti sono certamente soggetti a fluttuazioni, ci si deve domandare se le variazioni osservate rientrino entro limiti ragionevoli per un tasso fisso. Per poterlo decidere, occorre conoscere la probabilità di contare \(x\) eventi in un tempo fissato \(t\,\),
\begin{equation}\label{eq_2}
P_t(x)=\lim_{N_t\to\infty}\frac{N_t(x)}{N_t}\,,    
\end{equation} 
dove \(N_t(x)\) è il numero di prove ripetute in cui sono stati osservati esattamente \(x\) eventi nel tempo \(t\) su \(N_t\) prove totali (ciascuna di durata \(t\)). Vale evidentemente la condizione: 
\[\sum_{x=0}^\infty P_t(x)=1\,.\]
La funzione che esprime questa probabilità è, appunto, la poissoniana
\begin{equation}\label{eq_3}
P(x;\mu)=\frac{\mu\,^xe^{-\mu}}{x\,!},    
\end{equation}
che rappresenta la probabilità di conteggiare esattamente \(x\) eventi (sempre un numero intero) nel tempo \(t\) quando il numero atteso è \(\mu=\lambda\times t>0\) (generalmente non intero). È facile, infatti, verificare dall'eq. \ref{eq_3} che
\[\bar{x}=\sum_{x=0}^\infty x\times P(x;\mu)=\mu\,,\]
mentre la deviazione standard è semplicemente \(\sqrt{\mu}\,\), la radice quadrata della media. 

Per ragioni storiche, inerenti alla derivazione della poissoniana come limite della binomiale quando la probabilità di ciascun evento è piccola, la poissoniana è talvolta anche chiamata, a seconda della lingua degli autori, \textit{law of rare events}, \textit{loi des petites probabilités} o \textit{Gesetz der seltenen Ereignisse}. 
\section{Il processo di Poisson sulla linea temporale}\label{sez_3}
Le assunzioni che definiscono un processo di Poisson, limitandoci per ora al tipo più semplice, possono essere espresse mediante un piccolo numero di assiomi dai quali è possibile ricavare direttamente la poissoniana, come vedremo nella sezione successiva. Il sistema di postulati qui scelto (ripreso, con adattamenti, dal volume di Parzen \cite{parzen1962} e, in parte, da Ogborn e altri \cite{ogborn2003b}) comprende quattro assiomi: 
\begin{enumerate}
\item \textbf{Il n. di eventi \(x(t)\) evolve per incrementi indipendenti ("ipotesi di indipendenza")}. \\ Cioè, i conteggi di eventi in intervalli disgiunti sono variabili casuali indipendenti. Questa indipendenza consente, in ultima analisi, di calcolare la probabilità in un intervallo applicando le semplici regole di composizione delle probabilità di eventi indipendenti in sottointervalli.
\item \textbf{Per ogni \(t>0\), vale la condizione ("i. delle probabilità positive") \(0<P[x(t)>0]<1\,.\)}\\ Ossia, in parole, in ogni intervallo c'è una probabilità positiva che un evento possa accadere, ma non v'è certezza che accada. Il caso \(x=0\) (nessun accadimento) ha convenzionalmente probabilità 1 a \(t=0\,\) (\(P[x(0)=0]=1\,\)) e viene considerato talvolta come assioma a sé stante (assioma 0). Nell'intervallo semiaperto \(\interval({0,t}]\) vale anche la condizione {\(0<P[x(t)=0]<1\,\),} ovvero, c'è anche una probabilità positiva che non accada nessun evento in un tempo finito.
\item \textbf{Per ogni \(t\geq 0\,\), vale la condizione ("i. di rarità")
\(\lim\limits_{\tau\to 0}\frac{P[x(t+\tau)-x(t)\geq 2]}{P[x(t+\tau)-x(t)=1]}=0\,.\)}\\Detto a parole, in intervalli sufficientemente piccoli può ricadere, al più, un punto--evento; o anche, due o più eventi non possono accadere esattamente nello stesso istante. Intuitivamente, tornando alla fig. \ref{fig:PointProcess}, non vedremo mai i punti--evento "addossati" l'un l'altro, ma saranno separati.  Questo requisito di isolamento esige che, se si è verificato un evento in un intervallo sufficientemente piccolo, la probabilità condizionale di un secondo o di più eventi è trascurabile. Sono date solo due possibilità: o accade un evento o non accade nulla; la loro probabilità congiunta è 1, cioè, una certezza. Quindi, detta \(dp\) la probabilità che in un piccolo intervallo di tempo \(dt\) accada un evento, \(1-dp\)
 sarà la probabilità che non accada nulla (a meno di infinitesimi di ordine superiore in \(dt\)). Ciò esprime la condizione di "puntualità" di un processo poissoniano: gli eventi accadono "in un punto", intendendo con ciò un tempo molto minore di ogni intervallo di osservazione. 
\item \textbf{Il n. di eventi \(x(t)\) evolve per incrementi stazionari ("i. di equidistribuzione")}, \\ cioè, le variabili casuali che esprimono gli incrementi in intervalli equiestesi (non necessariamente disgiunti) sono equidistribuite.  Questo significa, in ultima analisi, che i punti--evento sono uniformemente distribuiti nell'intervallo di osservazione.
\end{enumerate}
L'ultimo assioma rinvia al concetto di distribuzione uniforme di eventi, che può essere meglio compreso introducendo la "densità di probabilità per evento": 

\textit{la probabilità infinitesimale \(dp\) che accada un evento in un intervallo infinitesimale \(dt\) è proporzionale al tasso \(\lambda\) (cioè, \(dp=\lambda \,dt\)), a meno di infinitesimi di ordine superiore in \(dt\). }

L'interpretazione di $\lambda$ che ne deriva, considerata intuitiva da taluni autori \cite{chiu2013}, è la seguente: suddividendo ciascuno degli intervalli \(t_i\,\), introdotti nella sez. \ref{sez_2}, in \(n_i\) intervallini equiestesi di grandezza \(\tau\) piccola a piacere, si può esprimere l'eq. \ref{eq_1} come
\[\lim\limits_{\substack{n\to\infty\\\tau\to 0}}\frac{x/n}{\tau}=\lambda\,,\]
dove il rapporto \(x/n\,\) al numeratore, in virtù dell'assioma 3, diventa 
\[\frac{x}{n}=\frac{(\text{intervallini in cui si è verificato un evento})}{(\text{intervallini totali})};\]
questo rapporto, passando al limite per \(n\to\infty,\,\tau\to 0\,\), si può interpretare come la probabilità infinitesimale \(dp\) dell'evento nell'intervallo infinitesimale \(dt\), e la corrispondente densità sarà data da
\begin{equation}\label{eq_4}
\frac{dp}{dt}=\lambda\,.    
\end{equation}
Il tasso \(\lambda\) assume nell'eq. \ref{eq_4} anche il significato di "intensità" (che in contesti più generali può riuscire funzione di $t$) e rappresenta la "probabilità per unità di tempo per evento" (coincidente con la densità di probabilità nel caso unidimensionale): ad es., è la probabilità che una particella che si muove a caso all'interno di un recipiente in cui sono praticati dei piccoli fori trovi una via di fuga in un tempo dato \cite{aloisi2018}. Dalla sua relazione con il valore atteso \(\mu=\lambda \,t\,\), introdotta nella sez. \ref{sez_2}, si vede che, quando \(\mu=1\), \(\lambda\) è l'inverso di un tempo caratteristico e il suo inverso \(\lambda^{-1}\), a sua volta, è il tempo tra due eventi consecutivi (tempo di interarrivo)\footnote{Per brevità viene omessa la trattazione dei tempi di interarrivo (o distanze tra eventi consecutivi) e della loro distribuzione.}, o il tempo di fuga della particella nell'esempio fatto. La costanza di \(\lambda\) rende la densità di probabilità indipendente dalla posizione dell'evento nell'intervallo di osservazione, come esige l'assioma 4, sicché sarà sempre possibile scegliere indipendentemente \(x\) punti \(t_1, \,t_2, \,\cdots, \,t_x\) in un intervallo \(\interval({0,t}]\) in modo che, in media, gli eventi siano uniformemente distribuiti su \(\interval({0,t}]\). Per semplicità (e per comodità) considereremo nel seguito l'eq. \ref{eq_4} e l'assioma 4 equivalenti, o intercambiabili. Supponendo che nel tempo \(t\) siano accaduti \(x\) eventi in \(x\) punti distinti della linea temporale, ciascuno con densità condizionale costante \(1/t\,\), si può dimostrare \cite{doob1990} che la densità condizionale di \(x\) punti-evento, scelti indipendentemente, di trovarsi in \(t_1, \,t_2, \,\cdots, \,t_x\) è \(x\,!/t\,^x\). Torneremo su questo risultato nella sez. \ref{sez_5}.

L'eq. \ref{eq_4} è talora assunta come punto di partenza per introdurre il processo di Poisson e viene usata nella derivazione diretta della poissoniana a partire dagli assiomi che lo definiscono. Ciò suggerisce che essa è una conseguenza inevitabile della "completa randomicità" insita nelle assunzioni (talvolta utilizzate in modo implicito) di indipendenza codificate in quegli assiomi \cite{kingman1993}. È interessante notare che facendo cadere uno o più degli assiomi 1--4 si possono caratterizzare processi stocastici più generali (ad es., rinunciando all'assioma 3 il processo risultante non ha una distribuzione poissoniana) o varianti del processo di Poisson in cui compare la funzione intensità \(\lambda(t)\) (ad es., sostituendo o modificando opportunamente l'assioma 4 è possibile modulare l'intensità nel tempo) \cite{haight1967}. Nel seguito consideriamo soltanto processi di Poisson con \(\lambda\) costante. Processi di questo tipo sono detti omogenei.
\section{Derivazione diretta della poissoniana}\label{sez_4}
Omettiamo la derivazione mediante il limite di Poisson, brevemente descritta all'inizio, trattandosi di un argomento assai noto che compare in molti libri di testo, sia di matematica che di fisica (v., ad es., \cite{reif1974}). Tralasciamo anche approcci più sofisticati, come quello della \textit{funzione generatrice} (anche questa ottenibile a partire dagli assiomi 1--4), o ancora quello della \textit{matrice di transizione} del processo di Poisson. Un metodo più consueto, che ugualmente omettiamo giacché adeguatamente descritto in molti libri e articoli (v., ad es., \cite{daboni1980}), fa uso delle equazioni differenziali--alle differenze del processo di Poisson: queste equazioni, una versione differenziale -- nota anche come master equation -- dell'equazione di Chapman–Kolmogorov, vengono generalmente risolte ricorsivamente. 

In questa sede presentiamo invece un metodo di derivazione diretta - o "esatta" - che prende le mosse dagli assiomi del processo di Poisson stazionario in una dimensione e fa uso soltanto di nozioni elementari di teoria delle probabilità. Nonostante la semplicità concettuale che lo contraddistingue, questo approccio non viene normalmente proposto nei libri di testo, con poche eccezioni\footnote{Fa eccezione, ad es., Yost \cite{yost1985}, che fornisce gli enunciati degli assiomi (una variante di quelli da noi scelti) ma omette la dimostrazione dettagliata.}. Una derivazione esatta si trova, in forma piuttosto concisa, nel volume di Bevington e Robinson \cite{bevington2003}, che a loro volta la traggono da Orear \cite{orear1982}. Si può vedere anche una discussione nel forum online \href{https://physics.stackexchange.com/q/372602}{\texttt{physics.stackexchange}}
  \cite{DanielSank}. 
  
  Facciamo a questo punto una breve digressione la cui utilità diverrà chiara più avanti: Il fattoriale \(x\,!\) nell'eq. \ref{eq_3} scaturisce, nel limite di Poisson, dal coefficiente binomiale \(\binom{n}{x}\,{\xRightarrow[\substack{}]{}}\,n^x/x\,!\), che a sua volta proviene, come noto, dall'enumerazione delle possibili combinazioni bernoulliane successo/insuccesso nel discreto (v., ad es., \cite{reif1974}). Nella derivazione esatta il fattoriale è invece il risultato di una integrazione nel continuo. Come vedremo meglio nella sez. \ref{sez_5}, quanto si interpretano fisicamente i risultati \(x\,!\) può anche essere interpretato
come il fattore combinatorio che esprime le possibili disposizioni di \(x\) particelle in altrettanti punti nell'intervallo \(\interval({0,t}]\). 

La derivazione diretta procede, in linea di massima, nel modo seguente: Partendo dalla nozione di "densità di probabilità per evento", introdotta nella sezione precedente, si suddivide l'intervallo \(\interval({0,t}]\) in \(n\) sottointervalli disgiunti di grandezza \(\tau=t/n\) piccola a piacere, indicati con \(dt_1, \,dt_2,\,\cdots, \,dt_n\,\). Si applica quindi l'eq. \ref{eq_4} (e, di riflesso, l'assioma 4) in combinazione con gli assiomi 3 ("rarità"), 1 ("indipendenza") e, implicitamente, con il 2 ("probabilità positive"),  per calcolare la densità di probabilità congiunta di \(x\) eventi direttamente come prodotto delle piccolissime probabilità elementari negli intervallini in cui è stato suddiviso \(\interval({0,t}]\).
\subsection{Probabilità di non osservare nessun evento nel tempo \textit{t}}
In ogni intervallino \(dt_i\) (\(i=1,\,2,\,\cdots,\,n)\) di \(\interval({0,t}]\) la probabilità infinitesima per evento è, in virtù dell'eq. \ref{eq_4}, \(\lambda \,dt_i\,\). L'evento, o accadrà o non accadrà e, per l'assioma 3 di "rarità", si può trascurare la probabilità che nell'intervallino possa accadere più di un evento. Pertanto, vi sarà una probabilità \(1-\lambda \,dt_i\) che l'evento non accada. Poiché, per l'assioma 1, gli eventi nei diversi intervallini sono indipendenti, la probabilità \(P(0,t,\lambda)\) che non accada nessun evento in nessun intervallino sarà quindi il \(\lim\limits_{n\to\infty}\) del prodotto
\[\prod_{i=1}^ n(1-\lambda \,dt_i) =(1-\frac{\lambda \,t}{n})^n=e^{n\ln(1-\lambda \,t/n)}\quad{\xRightarrow[(n\to\infty)]{}}\quad e^{-\lambda \,t},\]ovvero
\begin{equation}\label{eq_5}
P(0,t,\lambda)=e^{-\lambda \,t}.
\end{equation}
Si noti che \(P(0,0,\lambda)=1\) in accordo con l'assioma 2.
\subsection{Probabilità di un evento singolo}
Supponiamo che l'evento accada in un intervallino \(dt_i\) e non accada nient'altro nel periodo di osservazione \(t\). La probabilità elementare \(dP(1,t,\lambda)\) di un evento singolo è allora data dalla probabilità per evento \(\lambda \,dt_i\) nell'intervallino \(dt_i\) per la probabilità \(P(0,t,\lambda)\), espressa dall'eq. \ref{eq_5}, che non accada nessun altro evento nel tempo rimanente \(t-dt_i\) (che possiamo assumere uguale a \(t\) essendo \(dt_i\) infinitesimo),
\[dP(1,t,\lambda)=\lambda \,dt_i\times P(0,t,\lambda)\,,\]
e integrando su \(dt_i\) da 0 a \(t\,\),
\begin{equation}\label{eq_6}
P(1,t,\lambda)=\lambda \,t \,P(0,t,\lambda)=\lambda \,t \,e^{-\lambda \,t}.
\end{equation}
\subsection{Probabilità di eventi multipli} \label{ss_mult}
Consideriamo ora la distribuzione di \(x\) eventi in \(\interval({0,t}]\). Supponiamo che il primo evento accada a \(t_1\), il secondo a \(t_2\), ecc. e che \(0<t_1<t_2<\cdots<t_x\leq t\,\); si potrà allora considerare una serie di \(x\) intervallini disgiunti \(dt_1, \,dt_2,\,\cdots, \,dt_x\) entro i quali ricadono, rispettivamente, i punti--evento \(t_1, \,t_2,\,\cdots,\,t_x\). La probabilità elementare congiunta di questi \(x\) eventi \(d^xP(x,t,\lambda)\) è data dal prodotto delle probabilità \(\lambda \,dt_i\) per evento in ciascun intervallino della serie, per la probabilità \(P(0,t,\lambda)\) che non accada nessun altro evento nel tempo rimanente (avendo nuovamente assunto \(t-dt_1-dt_2-\cdots-dt_x\) uguale a \(t\)):
\begin{equation}\label{eq_7}
d^xP(x,t,\lambda)=\prod_{i=1}^ x\lambda \,dt_i\times P(0,t,\lambda)=\lambda^x\,e^{-\lambda \,t}\prod_{i=1}^ xdt_i\,.
\end{equation}La probabilità \(P(x,t,\lambda)\) di osservare \(x\) eventi nel periodo di osservazione \(t\) si ottiene integrando l'elementino di probabilità \(d^xP(x,t,\lambda)\) sul dominio dei possibili valori di \(t_1, \,t_2,\,\cdots,\,t_x\): 
\[P(x,t,\lambda)=\lambda^x\,e^{-\lambda \,t}\idotsint\limits_{0<t_1<\cdots<t_x\leq t} \,dt_1 \cdots \,dt_x\,.\]
Si tratta di un integrale definito, con condizioni nei limiti, della funzione \(f(t_1,\,\cdots,\,t_x)=1\,\), e rappresenta il volume di un \textit{x}--simplesso retto regolare (talvolta anche chiamato iperottante retto), la generalizzazione a \textit{x} dimensioni di un tetraedro trirettangolo in 3D (o di un triangolo rettangolo isoscele in 2D). Questo integrale di volume \textit{x}--dimensionale può essere calcolato scrivendolo sotto forma di integrale multi-nidificato (o multi-iterato) e integrando step-by-step partendo dall'integrale più interno,\[\idotsint\limits_{0<t_1<\cdots<t_x\leq t} \,dt_1 \cdots \,dt_x=\int\limits_0^t dt_1\int\limits_{t_1}^t dt_2\cdots\int\limits_{t_{x-1}}^t dt_x=\int\limits_0^t dt_1\int\limits_{t_1}^t dt_2\cdots\int\limits_{t_{x-2}}^t(t-t_{x-1}) dt_{x-1}=\cdots\]
(dove ci siamo fermati dopo il primo step).

Applicando ricorsivamente la regola delle potenze (e la regola della catena per i termini negativi) con il procedimento illustrato, ad es., nel forum \href{https://math.stackexchange.com/q/2553384}{\texttt{math.stackexchange}} \cite{IsaacBrowne}, si ottiene infine il seguente risultato (v. anche la nota \ref{footnote_5}): 
\[\idotsint\limits_{0<t_1<\cdots<t_x\leq t} \,dt_1 \cdots \,dt_x=\frac{t\,^x}{x\,!}\](omettiamo per brevità la dimostrazione completa). Perciò, la probabilità di osservare \(x\) eventi nel tempo \(t\) è data da
\begin{equation}\label{eq_8}
P(x,t,\lambda)=\frac{\lambda^x \,t\,^x}{x\,!}\,P(0,t,\lambda)=\frac{(\lambda\,t)\,^x\,e^{-\lambda \,t}}{x\,!},
\end{equation}o, in forma più stringata, dall'eq. \ref{eq_3} sostituendo nell'eq. \ref{eq_8} il valore atteso \(\mu=\lambda\,t\,\). 
\section{Distribuzioni condizionali di eventi e di particelle}\label{sez_5}
È immediato verificare, dalle eq.ni \ref{eq_7} e \ref{eq_8}, il risultato, cui si è accennato nella sez. \ref{sez_3}, sulla densità condizionale di \(x\) eventi in \(t_1, \,t_2, \,\cdots, \,t_x\) dati esattamente \(x\) eventi nell'intervallo \(\interval({0,t}]\,\),
\[p(t_1, \,t_2, \,\cdots, \,t_x\mid x,t)=\frac{d^xP(x,t,\lambda)/dt_1\cdots dt_x}{P(x,t,\lambda)}=\frac{\lambda^x\,P(0,t,\lambda)}{P(x,t,\lambda)}=\frac{x\,!}{t\,^x}\,,\]
quando ciascun evento ha densità condizionale costante \(1/t\) (come si vede immediatamente dalla formula precedente per il caso \(x=1\)). Che la distribuzione dei punti-evento sia uniforme sull'iperottante retto \(0<t_1<t_2<\cdots<t_x\leq t\,\) lo si vede in modo più intuitivo esprimendo il precedente risultato come
\[p(t_1, \,t_2, \,\cdots, \,t_x\mid x,t)=\frac{1}{\frac{t\,^x}{x\,!}}=\frac{1}{V_x^{oct}(t)}\,,\]
dove \(V_x^{oct}(t)\) è il "volume" (generalizzato a \textit{x} dimensioni) dell'iperottante: in 1D (linea) \(V_1^{oct}(t)=t\,\), in 2D (triangolo) \(V_2^{oct}(t)=\frac{t\,^2}{2}\,\), in 3D (tetraedro) \(V_3^{oct}(t)=\frac{t\,^3}{6}\,\), eccetera\footnote{L'interpretazione geometrica di questo risultato può essere resa più evidente esprimendo, in maniera alternativa all'eq. \ref{eq_7}, la probabilità congiunta mediante la formula ricorsiva\[P(x+1;\mu)=\int\limits_{0}^\mu P(x;\mu')\,d\mu'\,P(0;\mu-\mu')=\int\limits_{0}^\mu P(x;\mu')\,e^{-(\mu-\mu')}\,d\mu'=e^{-\mu}\int\limits_{0}^\mu P(x;\mu')\,e^{\mu'}\,d\mu'\]
(dove \(\mu=\lambda\, t\), \(\mu'=\lambda\, t'\) e \(t'\) è un qualche punto sulla linea temporale che si lascia a sinistra \(x\) punti-evento nell'intervallo \(\interval({0,t}]\)), che può essere ottenuta come soluzione analitica delle equazioni di Kolmogorov \cite{logan2009} o direttamente dagli assiomi 1--4. Definendo allora \(V_x^{oct}(\mu)\equiv P(x;\mu)\,e^{\mu}\) e tenendo conto che \(V_1^{oct}(\mu)=\mu\,\), l'integrale più a destra diventa
\[V_{x+1}^{oct}(\mu)=\int\limits_{0}^\mu V_{x}^{oct}(\mu')\,d\mu'\,,\]
di immediato significato geometrico, da cui è facile ricavare per iterazione 
\[V_{x}^{oct}(\mu)=\frac{[V_1^{oct}(\mu)]^x}{x\,!}=\frac{\mu^x}{x\,!}\,.\]
}\label{footnote_5}. L'unità al numeratore rappresenta la probabilità condizionale di \(x\) eventi in punti \(t_1, \,t_2, \,\cdots, \,t_x\) dati \(x\) eventi nell'intervallo \(\interval({0,t}]\). Secondo un'interpretazione di taluni autori \cite{daley2003}, apparentemente fatta propria anche da \cite{bevington2003} e \cite{orear1982}, nell'ottenere questo risultato i punti sono trattati come indistinguibili eccetto per le loro posizioni. Ma in situazioni fisiche questi punti potrebbero essere le posizioni di \(x\) particelle fisicamente distinguibili e il fattoriale \(x\,!\) che origina in prima istanza come (inverso del) volume dell'iperottante unitario (\([V_x^{oct}(1)]^{-1}\)) può essere anche interpretato come il fattore combinatorio che rappresenta il numero di modi in cui \(x\) particelle distinte possono essere posizionate in \(x\) punti distinti dell'intervallo considerato. Le particelle individuali devono essere allora immaginate come uniformemente e indipendentemente distribuite su \(\interval({0,t}]\). Seguendo quest'interpretazione, la funzione che rappresenta la distribuzione di \(x\) particelle distinte su questo intervallo verrebbe scritta come
\[p(t_1, \,t_2, \,\cdots, \,t_x\mid x,t)=\frac{x\,!}{V_x^{rec}(t)}\,\]
dove ora \(V_x^{rec}(t)\) è il volume dell'iperrettangolo \((0<t_1\leq t\,,\cdots\,,0<t_x\leq t)\,\), ovvero \(t\,^x\). In questo senso, si  può dire che le densità condizionali poissoniane di \(x\) eventi uniformemente distribuiti sull'iperottante \(0<t_1<t_2<\cdots<t_x\leq t\,\) corrispondono alle densità di probabilità di \(x\) particelle distinte uniformemente distribuite, a caso, sull'intervallo \(\interval({0,t}]\). Incidentalmente, osserviamo che sfruttando la relazione, puramente geometrica \cite{bengtsson2017}, tra i volumi dell'iperottante e dell'iperrettangolo, 
\[V_x^{oct}=\frac{1}{x\,!}V_x^{rec}\,,\]
è possibile evitare il calcolo esplicito (...un po' tedioso!) dell'integrale multiplo sull'iperottante, esprimendolo tramite un integrale separabile (sull'iperrettangolo) che può essere fattorizzato con un prodotto di integrandi indipendenti. Infatti, 
\[\textcolor{red}{\idotsint\limits_{\substack{0<t_1<\cdots<t_x\leq t\\(x-\texttt{ottante})}}\,dt_1 \cdots \,dt_x}=\frac{1}{x\,!}\idotsint\limits_{\substack{0<t_1\leq t\\\cdots\\0<t_x\leq t\\(x-\texttt{rettangolo})}}\,dt_1 \cdots \,dt_x=\frac{1}{x\,!}\int\limits_0^t dt_1\int\limits_0^t dt_2\,\cdots\int\limits_0^t dt_x=\textcolor{red}{\frac{1}{x\,!}\prod_{i=1}^ x\int\limits_0^t dt_i}\]
e, pertanto,
\[P(x,t,\lambda)=\lambda^x\,e^{-\lambda \,t}\times\textcolor{red}{\idotsint\limits_{0<t_1<\cdots<t_x\leq t} \,dt_1 \cdots \,dt_x}=\lambda^x\,e^{-\lambda \,t}\times\textcolor{red}{\frac{1}{x\,!}\prod_{i=1}^ x\int\limits_0^t dt_i=}\]\[=\lambda^x\,e^{-\lambda \,t}\times\frac{1}{x\,!}\biggl[\int\limits_0^t dt\biggr]^x=\lambda^x\,e^{-\lambda \,t}\times\frac{1}{x\,!}t\,^x\,,\]
da cui, infine, l'eq. \ref{eq_8}.
\section{Distribuzione binomiale delle particelle in un sottointervallo}\label{sez_6}
Dalle considerazioni precedenti, o direttamente dall'eq. \ref{eq_8}, non è difficile verificare che il numero delle particelle in un sottointervallo di \(\interval({0,t}]\) segue una distribuzione binomiale quando il numero totale delle particelle nell'intervallo \(\interval({0,t}]\) è fissato. Infatti, suddividendo \(\interval({0,t}]\) in \(k\) sottointervalli, disgiunti tranne forse per alcuni estremi, \(\Delta t_1\,,\Delta t_2\,,\cdots\;,\Delta t_k\) (non necessariamente piccoli) tali che \(\Delta t_1+\cdots+\Delta t_k\leq t\,\), allora la probabilità congiunta di trovare \(x_1\,,\ldots\,,x_k\) particelle uniformemente distribuite, rispettivamente, in \(\Delta t_1\,,\ldots\,,\Delta t_k\,\) e trovarne \(\,x_r=n-x_1-\cdots-x_k\,\) nella parte rimanente \(t-\Delta t_1-\cdots-\Delta t_k\,\) è data da
\[P(x_1,\,\cdots,\,x_k\,;\,x_r)=\prod_{i=1}^ k\lambda\,^{x_i}\,e^{-\lambda\,\Delta t_i}\frac{1}{x_i\,!}(\Delta t_i)\,^{x_i}\times\]
\[\times\lambda\,^{x_r}\,e^{-\lambda\,(t-\Delta t_1-\cdots-\Delta t_k)}\frac{1}{x_r\,!}(t-\Delta t_1-\cdots-\Delta t_k)\,^{x_r}\,=\]
\[=\lambda\,^{n}\,e^{-\lambda\,t}\frac{1}{{x_1}\,!\cdots\,{x_k}\,!\times x_r\,!}(\Delta t_1)\,^{x_1}\cdots(\Delta t_k)\,^{x_k}(t-\Delta t_1-\cdots-\Delta t_k)\,^{x_r}\,,\]
e la corrispondente probabilità condizionale, fissato il numero totale \(n\) di particelle, da
\[P(x_1,\,\cdots,\,x_k\,;\,x_r\mid n)=\frac{P(x_1,\,\cdots,\,x_k\,;\,x_r)}{P(n,t,\lambda)}=\]
\[=\frac{\cancel{\lambda\,^{n}\,e^{-\lambda\,t}}\times\frac{1}{{x_1}\,!\cdots\,{x_k}\,!\,\times\, x_r\,!}(\Delta t_1)\,^{x_1}\cdots(\Delta t_k)\,^{x_k}(t-\Delta t_1-\cdots-\Delta t_k)\,^{x_r}}{\cancel{\lambda\,^{n}\,e^{-\lambda\,t}}\times\frac{1}{n\,!}t\,^n}=\]
\[=\frac{n\,!}{{x_1}\,!\cdots\,{x_k}\,!\times x_r\,!}(\frac{\Delta t_1}{t})\,^{x_1}\cdots(\frac{\Delta t_k}{t})\,^{x_k}(1-\frac{\Delta t_1}{t}-\cdots-\frac{\Delta t_k}{t})\,^{x_r}\,,\]
dopo aver semplificato i fattori comuni e diviso numeratore e denominatore per \(t\,^n\). Ora, interpretando le \(p_i=\Delta t_i/t\,\), per \(i=1\,,\cdots\,,k\,,\) come le probabilità di osservare le particelle individuali (che, come abbiamo visto nella sezione precedente, devono essere pensate uniformemente distribuite, a caso, su \(\interval({0,t}]\)), in ciascuno dei sottointervalli \(\Delta t_i\,\) e la \(p_r=1-\Delta t_1/t-\cdots-\Delta t_k/t\,\) come la probabilità di osservarle nella parte rimanente, si riconosce nell'espressione risultante la funzione di probabilità della distribuzione multinomiale di parametri \(((p_1,\,\cdots,\,p_r),\,n)\) con \(p_1+\cdots+p_r=1\,\),
\[P(x_1,\,\cdots,\,x_r\mid n)=M(x_1,\,\cdots,\,x_r;\,n,\,p_1,\,\cdots,\,p_r)=\binom{n}{x_1,\dots,x_r}\prod_{i=1}^ r p_i\,^{x_i}\,,\]
dove il coefficiente multinomiale
\[\binom{n}{x_1,\,\dots,\,x_r}=\frac{n\,!}{{x_1}\,!\,\cdots\,{x_r}\,!}\]
generalizza il coefficiente binomiale al numero di possibili sequenze con \(x_1,\,\ldots,\,x_r\) particelle, rispettivamente, in \(\Delta t_1,\,\ldots,\,t-\Delta t_1-\cdots-\Delta t_k\,\).
Nel caso si consideri un unico sottointervallo \(\Delta t\) ci si riduce alla funzione 

\[P(x,\,n-x\mid n)=B(x,n)=\binom{n}{x}p\,^x q\,^{n-x}\,,\]
con \(p=\Delta t/t\) e \(q=1-p=1-\Delta t/t\,\), che rappresenta la funzione di probabilità della distribuzione binomiale di \(x\) di particelle nel sottointervallo \(\Delta t\) fissatone il numero \(n\) nel più ampio intervallo \(\interval({0,t}]\).

Questa equivalenza tra un processo binomiale e un processo poissoniano omogeneo, condizionato sul numero di eventi, è una conseguenza di un teorema generale di equivalenza tra processi puntuali \cite{schabenberger2017}.
\section{La "riscoperta" della poissoniana}
È interessante notare, rinviando, ad es., al testo di Last e Penrose \cite{last2017} per una rassegna storica più ampia, come agli inizi del Novecento sia avvenuta una serie di "riscoperte" indipendenti della poissoniana nel contesto di applicazioni pratiche. All'astronomo svedese Carl V. L. Charlier, pioniere della statistica astronomica, viene accreditata la prima derivazione dell'eq. \ref{eq_8} mediante le equazioni differenziali--alle differenze (master equation) del processo di Poisson \cite{charlier1905}. A destare maggiore attenzione fu però la derivazione fornita da H. Bateman, in appendice a un articolo di Rutherford e Geiger sulle fluttuazioni nel conteggio delle particelle \(\alpha\) emesse da sostanze radioattive \cite{rutherford1910}, in quanto direttamente legata a un processo fisico fondamentale. Bateman mostrò che il numero di particelle \(\alpha\) emesse in periodi di tempo fissati soddisfaceva a un semplice insieme di equazioni differenziali le cui soluzioni erano le probabilità di Poisson \cite{bateman1910}. Il metodo di Bateman ha avuto largo seguito in sviluppi successivi (v., ad es., \cite{hughes1939}), ma è bene precisare che l'applicazione della statistica poissoniana nella descrizione del decadimento radioattivo è solo un'approssimazione, che può funzionare più o meno bene a seconda delle situazioni sperimentali \cite{sitek2015}. Passando ad altri campi della scienza, troviamo applicazioni molto note, come la classica analisi di Erlang del traffico telefonico \cite{erlang1909}, accanto ad altre meno note, come le ricerche in campo sanitario di A. G. McKendrick, ufficiale medico dell'armata britannica in India. Nei suoi primi due articoli compare la derivazione della poissoniana mediante la master equation e, a testimonianza della versatilità di questo metodo (...e di questo studioso!), l'applicazione al conteggio 
delle cellule nella fagocitosi leucocitaria \cite{mckendrick1914a,mckendrick1914b}, alle collisioni molecolari di un gas \cite{mckendrick1914a} e alle reazioni chimiche monomolecolari \cite{mckendrick1914a}. Gli esempi potrebbero essere numerosi, ma il punto che vogliamo evidenziare è che, indipendentemente dal particolare trattamento matematico, l'approccio assiomatico (insieme con l'eq. \ref{eq_4}) che abbiamo discusso nella sez. \ref{sez_3} si rivela particolarmente fecondo nelle applicazioni modellabili con processi poissoniani: se, nella situazione concreta, possiamo essere sicuri che ricorrono gli assiomi 1--4 (o loro controparti spazio-temporali), allora ne segue necessariamente l'eq. \ref{eq_8} (o la \ref{eq_3}) e la statistica è poissoniana. Vedremo una esemplificazione di questo metodo nella prossima sezione.
\section{Applicazione al gas perfetto classico}
Come esempio di applicazione dell'approccio assiomatico consideriamo il seguente problema: 

\textit{dimostrare che la probabilità di trovare \(n\) particelle in una data regione di volume \(v\) all'interno di un volume più ampio \(V\) (un "serbatoio") di un gas perfetto a temperatura \(T\) è data dalla poissoniana 
\[p_n=\frac{\bar{n}^ne^{-\bar{n}}}{n\,!}\]
dove \(\bar{n}\) è il numero medio di particelle nel volume \(v\). }

Dato che il numero di particelle nel volume interno è variabile per effetto degli scambi (in equilibrio) con il serbatoio a temperatura costante, il problema può essere risolto con i metodi della meccanica statistica tramite la funzione di partizione gran canonica \cite{mandl1971}. In alternativa, si può dimostrare che il conteggio del numero di particelle nel volume dato è modellabile con un processo poissoniano, ovvero che gli assiomi 1--4 (o, meglio, una loro controparte tridimensionale) sono applicabili al sistema in esame. 
\subsection{Il processo di Poisson in 3 dimensioni}
Nel caso di un processo di Poisson omogeneo in 3 dimensioni, con \(n\) in luogo di \(x\) e \(v\) in luogo di \(t\,\), in cui il numero di particelle evolve man mano che il conteggio viene esteso a regioni sempre più ampie, gli assiomi da noi scelti potrebbero assumere la forma seguente:
\begin{enumerate}[i]
\item \textbf{Il n. di particelle \(n(v)\) evolve per incrementi indipendenti ("ipotesi di indipendenza")}. 
\item \textbf{Per ogni \(v>0\), vale la condizione "i. delle probabilità positive") \(0<P[n(v)>0]<1\)}. \\ Ossia, in parole, in ogni regione c'è una probabilità positiva, ma non la certezza, di trovare almeno una particella. Il caso \(n=0\) (regione vuota) ha convenzionalmente probabilità 1 per un volume nullo: \(P[n(0)=0]=1\,\)(assioma 0). Nell'intervallo semiaperto \(\interval({0,v}]\) vale, inoltre, la condizione \(0<P[n(v)=0]<1\,\), ovvero, c'è anche una probabilità positiva che una regione finita sia vuota.
\item \textbf{Per ogni \(v\geq 0\,\), vale la condizione ("i. di rarità") \(\lim\limits_{\nu\to 0}\frac{P[n(v+\nu)-n(v)\geq 2]}{P[n(v+\nu)-n(v)=1]}=0\,\)}. \\ Cioè, in regioni sufficientemente piccole si può trovare, al più, una particella. In altre parole, non è possibile che più particelle si trovino esattamente nello stesso punto: ogni particella avrà intorno a sé dello spazio vuoto. Questa proprietà di isolamento garantisce che se in una regione sufficientemente piccola si è trovata una particella, la probabilità condizionale di trovare una seconda o più particelle è trascurabile. 
\item \textbf{il n. di particelle \(n(v)\) evolve per incrementi stazionari ("i. di equidistribuzione")}, che significa, in ultima analisi, che ciascuna particella di gas è uniformemente distribuita in \(V\) indipendentemente dalle restanti particelle.
\end{enumerate}
Per verificare l'applicabilità degli assiomi i--iv è cruciale la definizione di gas perfetto: è essenziale che le molecole di gas si muovano pressoché liberamente nello spazio (assioma i), che siano abbastanza separate le une dalle altre (assioma iii), che interagiscano tra loro molto debolmente\footnote{Più precisamente, che l'energia d'interazione sia trascurabile rispetto all'energia cinetica.} (assioma iv).  
L'eq. \ref{eq_4} nella sez. \ref{sez_3} ha la propria controparte tridimensionale in
\begin{equation}\label{eq_9}
\frac{dp}{dv}=\lambda\,.    
\end{equation}
\subsection{Fluttuazioni di Poisson}
La distribuzione tridimensionale delle particelle si può a questo punto ricavare ripercorrendo, \textit{mutatis mutandis}, tutti passaggi della derivazione diretta della poissoniana fatta nella sez. \ref{sez_4}, ma con una avvertenza: a differenza degli eventi sulla linea temporale, non esiste un ordinamento naturale per le particelle nello spazio tridimensionale. Possiamo però immaginare la regione di interesse come un "cilindro" di sezione variabile \(s(t)\) e volume \(v=\int_0^t s(t)\,dt\,\), suddiviso in "dischi" paralleli di spessore infinitesimo, tali che la probabilità che in un disco si trovi più di una particella (supposta puntiforme) sia trascurabile (come esige l'assioma iii). Ciascun disco avrà un volume infinitesimo \(dv=s(t)\,dt\,\). Supponiamo, inoltre, che le particelle si trovino in \(n\) di questi dischi, posizionati lungo l'"asse" del cilindro nei punti \(t_1,\ldots,t_n\,\) con \(0<t_1<\cdots<t_n\leq t\,\). Avremo, allora, una relazione tra le posizioni delle particelle corrispondente all'ordinamento dei punti--evento sulla linea temporale e potremo integrare come nella sez. \ref{sez_5} (la dimostrazione è omessa)\footnote{In alternativa si può utilizzare l'integrale della nota \ref{footnote_5}, nella forma
\[p_{n+1}(\mu)=e^{-\mu}\int\limits_0^\mu p_n(\mu')\,e^{\mu'}\,d\mu'\]
dove ora \(\mu=\lambda\, v\), \(\mu'=\lambda\, v'\) e \(v'\) è (il volume di) una qualche sottoregione che racchiude \(n\) particelle nella più ampia regione \(v\). Avremo allora \(V_n^{oct}(\mu)\equiv p_n(\mu)\,e^\mu\,\), \[V_{n+1}^{oct}(\mu)=\int\limits_0^\mu V_n^{oct}(\mu')\,d\mu'\,,\] e infine \[V_n^{oct}(\mu)=\frac{\mu^n}{n\,!}\,.\]}:
\[\textcolor{red}{\idotsint\limits_{\substack{0<t_1<\cdots<t_n\leq t\\(n-\texttt{ottante})}}\,s(t)\,dt_1 \cdots \,s(t)\,dt_n}=\frac{1}{n\,!}\idotsint\limits_{\substack{0<t_1\leq t\\\cdots\\0<t_n\leq t\\(n-\texttt{rettangolo})}}\,s(t)\,dt_1 \cdots \,s(t)\,dt_n=\]
\[=\frac{1}{n\,!}\int\limits_0^t s(t)\,dt_1\int\limits_0^t s(t)\,dt_2\,\cdots\int\limits_0^t s(t)\,dt_n=\textcolor{red}{\frac{1}{n\,!}\prod_{i=1}^ n\int\limits_0^t s(t)\,dt_i}=\frac{1}{n\,!}\biggl[\int\limits_0^t s(t)\,dt\biggr]^n=\frac{1}{n\,!}v\,^n\,.\]
Otterremo, infine, come controparte dell'eq. \ref{eq_8}, la seguente:
\begin{equation}\label{eq_10}
p_n=\frac{(\lambda\,v)\,^n\, e^{-\lambda \,v}}{n\,!}\,.
\end{equation}
\subsection{Limite di Poisson in meccanica statistica}
In modo equivalente, poiché il numero \(N\) di molecole del volume grande \(V\,\) è fisso e le particelle sono individualmente distinguibili (caso classico), si possono applicare i risultati della sez. \ref{sez_6} sulla distribuzione binomiale delle particelle in un sottointervallo e scrivere
\[P(n,\,N-n\mid N)=B(n,N)=\binom{N}{n}p\,^n q\,^{N-n}\,,\]
con \(p=v/V\) e \(q=1-p=1-v/V\,\), che rappresenta la funzione di probabilità della distribuzione binomiale del numero \(n\) di particelle nel volume \(v\) fissatone il numero \(N\) nel più ampio volume \(V\). Come usuale in meccanica statistica, si assume che \(V\) e \(N\,\) tendano all'infinito in modo tale che la densità numerica delle particelle rimanga finita, \(\frac{N}{V}\,{\xRightarrow[\substack{}]{}}\,\lambda>0\, \) per \({\,V\to\infty\,,\,N\to\infty}\). Questa condizione viene talvolta chiamata limite termodinamico. Un'applicazione rigorosa del limite di Poisson in meccanica statistica, con una generalizzazione a dimensioni superiori, si può trovare nel bel volume di Sinai \cite{sinai1992}. Si può verificare abbastanza agevolmente (ne omettiamo la dimostrazione dettagliata) che eseguendo il passaggio al limite si ottiene
\[p_n=\binom{N}{n}(\frac{v}{V})\,^n (1-\frac{v}{V})\,^{N-n}\,{\xRightarrow[(V,N\to\infty)]{}}\,\frac{1}{n\,!}(\frac{Nv}{V})\,^ne^{-\frac{Nv}{V}}\,{\xRightarrow[(V,N\to\infty)]{}}\,\frac{1}{n\,!}(\lambda v)\,^ne^{-\lambda v}\,,\]
cioè, l'eq. \ref{eq_10}, che è la poissoniana con parametro \(\bar{n}=\lambda v\,\).

I risultati delle due ultime sottosezioni sono generalizzabili a regioni con un numero arbitrario di dimensioni spaziali \cite{sinai1992}.
\section{Considerazioni conclusive}
Il metodo di derivazione della poissoniana dal processo stocastico sottostante proposto in questo articolo si caratterizza, rispetto ad altri approcci indipendenti, per una maggiore semplicità, non richiedendo l'uso di equazioni differenziali e la loro soluzione, né presupponendo altre conoscenze avanzate di matematica e di statistica. È  sufficiente l'applicazione diretta di alcuni risultati elementari di teoria delle probabilità e di un piccolo numero di assiomi, il che rende il metodo applicabile a un'ampia varietà di problemi. In conclusione, da questo approccio il docente potrà trarre qualche spunto per una presentazione introduttiva della teoria delle probabilità a un livello proponibile, con le dovute cautele, in una classe liceale.    


\end{document}